%% file: modelarxiv.tex
\documentclass[12pt]{article}

  \usepackage{amssymb}         
  \usepackage{amsmath}          
  \usepackage{amsfonts}           
  \usepackage{amsthm}
 
  \usepackage[mathscr]{eucal}

  \usepackage{graphicx}        
  \usepackage{color}

  \newcommand{\per}{\operatorname{Per}_1}
  \newcommand{\inter}{\operatorname{int}}
  \newcommand*{\R}{\mathbb{R}}
  
  \newcommand*{\C}{\mathbb{C}}
  \newcommand*{\N}{\mathbb{N}}

  \newcommand*{\D}{\mathbb{D}}
  
  \newcommand*{\CC}{\hat\C}

\newtheorem{definition}{Definition}
\newtheorem{etheorem}{Theorem}

\newtheorem{lemma}{Lemma}

\title{A model for the parabolic slices $\per(e^{2\pi i p/q})$ in moduli
  space of quadratic rational maps}
\author{Eva Uhre}

\begin{document}

\maketitle
\begin{abstract}
The notion of \emph{relatedness loci} in the parabolic slices
  $\per(e^{2\pi i p/q})$ in moduli space of quadratic rational maps is introduced.
  They are counterparts of the disconnectedness or escape locus in the
  slice of quadratic polynomials. 
  A model for these loci is presented,
  and a strategy of proof of the faithfulness of the model is given.
\end{abstract}

\section{Introduction}
Let $\mathcal{M}_2$ denote the moduli space of M\"obius conjugacy
classes of quadratic rational maps $f\colon\CC\to\CC$. Following
definitions and statements from Milnor\mbox{ \cite{M93}}, consider loci:
\[\per(\lambda)=\left\{[f]\in \mathcal{M}_2 : f \mbox{ has a fixed point with
eigenvalue } \lambda\right\}\cong \C. \]

Here the focus will be on parabolic slices $\per(\omega)$, with
$\omega=e^{2\pi i p/q}$, $p/q\neq 0/1$, i.e. those consisting of equivalence
classes of maps with a parabolic fixed point with eigenvalue $\omega$.
In such a slice the dynamics is characterized according to the behavior of the critical
points. The \emph{relatedness} locus
$\mathcal{R}^{\omega}$ in $\per(\omega)$ is defined by:
\begin{equation} 
  \mathcal{R}^{\omega}=\left\{[f]\in \per(\omega):
  \lim_{n\to\infty}f^n(c_1)=z_0=\lim_{n\to\infty}f^n(c_2)\right\}, 
\end{equation} 
where $z_0$ is the (persistent) parabolic fixed point and $c_1$ and
$c_2$ are the critical points of $f$. The locus
$\mathcal{R}^{\omega}$ is neither open nor closed. It consists
of open, connected components of maps where both critical points are
in the parabolic basin (in \mbox{ \cite{M93}} called
\emph{bitransitive} and \emph{capture} components respectively,
according to whether both or only one critical point is in the
immediate basin), a countable set of points corresponding to maps
where one critical point is eventually mapped to the parabolic fixed
point, and a finite set of points corresponding to maps where the
parabolic fixed point is degenerate, i.e. has two $q$-cycles of
components in the immediate basin. In the slice $\per(0)$ the relatedness locus $\mathcal{R}^{0}$
is the escape locus $\C\setminus M$, where $M$ is
the Mandelbrot set, the connectedness locus in the slice of polynomials.

\section{The model}
The objective is to construct a model for
$\mathcal{R}^{\omega}$(see theorem \ref{thm:map}). Consider the
quadratic polynomial
\[P_{\omega}(z)=\omega z+ z^2,\] with a parabolic fixed point with
multiplier $\omega$ at 0. This fixed point is called the
\emph{$\alpha$-fixed point}. Let $\Lambda_{\omega}$ denote the parabolic
basin of 0 for $P_{\omega}$ and define also an augmented basin $\tilde
\Lambda_\omega$:  
\[\tilde \Lambda_{\omega}=\left\{z\in \CC:
  \lim_{n\to\infty}P_{\omega}^n(z)=0\right\} = \Lambda_{\omega}\cup
\left\{z\in \CC: \exists n\geq 0, P_{\omega}^n(z)=0\right\}.\] 
The immediate basin has $q$ components, labelled $B_j$, $j\in\{0,...,q-1\}$ counter--clockwise, so that $B_0$ contains the critical point $-\omega/2$.
It follows from the theory of quadratic polynomials that there are $q$
external rays landing at 0, dividing $\CC$ into $q$ components. Let
$S_p$ denote the component containing the critical value
$P_{\omega}(-\omega/2)$. Let $\phi_{\omega}\colon \Lambda_{\omega}\to\C$ be an extended Fatou
coordinate for $P_\omega^q$, i.e. a surjective holomorphic map, of
infinite degree, with critical points at the critical point
$-\omega/2$ of $P_\omega$ and at all its pre--images, so that
$\phi_{\omega}\circ P_{\omega}^q=1+\phi_{\omega}$.

Normalize $\phi_{\omega}$ so that $\phi_{\omega}\circ
P_{\omega}=1/q+\phi_{\omega}$ and $\phi_{\omega}(-\omega/2)=0$. Let
$\mathcal{P}^j_{\delta}\subset B_j$ be
the connected component of $\phi^{-1}_{\omega}(\{z=x+iy:
x>\delta\})\cap B_j$ with the fixed point 0 on its boundary, called a \emph{petal}. 

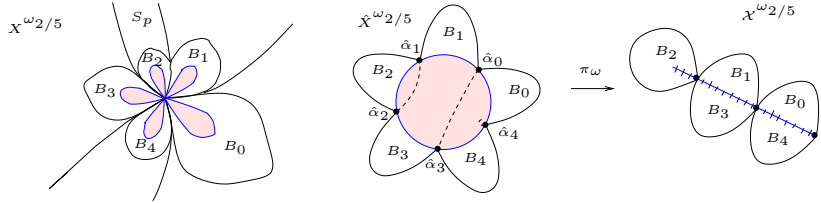
\begin{figure}[htbp]
  \centering
  \input{figure1.pstex_t}
\begin{quote}
  \caption{\label{fig:carcomp} A sketch of the construction of
    $\mathcal{X}^{\omega}$, in the case $p/q=2/5$.}
\end{quote}
\end{figure}

Denote by $X^{\omega}$ the subset of the Riemann sphere obtained by
removing the union  of the closures of the petals $\mathcal{P}_{1/q}$,
$X^{\omega}=\CC\setminus
\bigcup_{j=0}^{q-1}\overline{\mathcal{P}^j_{1/q}}$. Let $\hat
X^{\omega}\cong \CC\setminus\D$ be the Carath\'eodory
  compactification of $X^{\omega}$, i.e. the disjoint union of
$X^{\omega}$ and the set consisting of all prime ends of $X^{\omega}$.
The boundary of $\hat X^{\omega}$ can be naturally identified with the
boundaries of the petals, together with $q$ copies of the
$\alpha$-fixed point, corresponding to the $q$ different accesses to
the $\alpha$-fixed point from $X^{\omega}$. The copies are labelled $\hat \alpha_j$, $j\in\{0,...,q-1\}$ counter--clockwise, so that $\hat\alpha_j$ is an endpoint of $\partial
\mathcal{P}^j$ and $\partial \mathcal{P}^{j+1}$ (figure
\ref{fig:carcomp}). For the remainder of this section sets in $\CC\setminus \bigcup_{j=0}^{q-1}
\mathcal{P}^j_{1/q}$ are to be understood as subsets of $\hat
X^{\omega}$, i.e. with $q$ copies of the $\alpha$-fixed point. 

Now an equivalence relation on $\hat X^{\omega}$ is defined. To shorten notation let $\partial \mathcal{P}^j
= \partial \mathcal{P}^j_{1/q}$.
\begin{definition}
  Two points $z_1 \in \partial \mathcal{P}^j$ and $z_2 \in
  \partial \mathcal{P}^k$ are called \emph{equivalent modulo $p/q$}, written
  $z_1 \sim_{p/q} z_2$, if the following two conditions are
  satisfied:
\begin{itemize}
\item $j+k=2p \mod q$, and
\item $\phi_{\omega}(z_1)+\phi_{\omega}(z_2)=2/q$.
\end{itemize}
Two points $\hat\alpha_j$ and $\hat\alpha_k$ are said to be equivalent
modulo $p/q$ if $j+k=(2p-1) \mod q$.
\end{definition}

Let $\mathcal{X}^{\omega}$ be the quotient of $\hat X^{\omega}$ under
the equivalence relation $\sim_{p/q}$, let $\pi_\omega\colon \hat
X^{\omega}\to \mathcal{X}^{\omega}$ denote the projection map induced
by $\sim_{p/q}$ and let
$\nu_\omega=\pi_\omega(\partial \hat X^{\omega})\subset \mathcal{X}^{\omega}$ denote the scar
after gluing the real--analytic
boundaries of the petals back together under $\sim_{p/q}$. It can be proved\mbox{ \cite{U09}} that the map $\pi_\omega$ gives $\mathcal{X}^{\omega}$ a
  Riemann surface structure which extends the initial structure of
  $X^\omega$, and so that $\mathcal{X}^{\omega}\cong \CC$.

\begin{definition}
  The model space $\hat \Lambda^\omega\subset \mathcal{X}^{\omega}$
  for $\mathcal{R}^{\omega}$ is defined by:
  \[\hat
  \Lambda^{\omega}=\left (\tilde \Lambda_{\omega}\setminus (S_p \cup
    \bigcup_{j=0}^{q-1} \mathcal{P}^j_{1/q})\right)/_{\sim_{p/q}}.\]
\end{definition}

\begin{definition}
From the Fatou coordinate define a tree in $\hat
\Lambda^\omega$, called a \emph{bubble--tree} and denoted $\hat{\mathcal{T}}^\omega$, by:
\[\hat{\mathcal{T}}^\omega=\pi_\omega\left(\phi_\omega^{-1}(\R)
  \cup \bigcup_{n> 0}P_\omega^{-n}(0)\right) \cup\nu_\omega.\]
The bubble--tree has vertices at pre-fixed and (pre)-critical points of $P_\omega$ and at the points $\pi_\omega(\hat \alpha_i)$, $i\in\{0,...,q-1\}$. A metric is defined on the tree by assigning length one to every edge and letting the distance between any two vertices be the sum of the lengths of the edges in the unique finite path between them.  
\end{definition}

\section{Faithfulness of the model}
\begin{etheorem}\label{thm:map}
  There exists a bijective map $\chi^{\omega}\colon \mathcal{R}^{\omega}\to \hat \Lambda^{\omega}$, which is conformal in $\inter(\mathcal{R}^{\omega})$. The inverse $(\chi^{\omega})^{-1}$ is continuous on compact subsets of the bubble--tree
  $\hat{\mathcal{T}}^\omega$, with respect to the topology induced by the metric on the tree.
\end{etheorem}
Let $f_\sigma \in \mathcal{R}^\omega$ be non-degenerate parabolic, let $\phi_\sigma$ be a
Fatou coordinate for $f_\sigma$ and let $R\in\R$ be smallest so that
the union of petals $\bigcup_{j=0}^{q-1}\mathcal{P}^j_{\sigma,R}$ contains no
critical point, but contains (at least) one critical point on the
boundary. These petals are called maximal attracting petals and (one of) the critical point(s) on the boundary is called the closest critical point and denoted $c_1$. The
other critical point is then called the second critical point and denoted $c_2$. The critical values under $f_\sigma$ are denoted $v_1$ and $v_2$ respectively. Normalize the Fatou coordinate $\phi_\sigma$ so that $\phi_\sigma(c_1)=0$ and $\phi_{\sigma}\circ
f_{\sigma}=1/q+\phi_{\sigma}$. Let $\mathcal{U}_\omega^0 =
\bigcup_{j=0}^{q-1}\mathcal{P}_0^j$ be the maximal attracting petals
for $P_\omega$, and $U_\sigma^0=
\bigcup_{j=0}^{q-1}\mathcal{P}_{\sigma,0}^j$ the maximal attracting
petals for $f_\sigma$. Further, let
$U_\sigma^n=f_\sigma^{-n}(U_\sigma^0)$ and
$\mathcal{U}_\omega^n=P_\omega^{-n}(\mathcal{U}_\omega^0)$. The map $\chi^{\omega}$ is constructed via a dynamical conjugacy:
\begin{lemma}\label{lemma:conjugacy}
  For all non-degenerate parabolic $f_\sigma\in
  \mathcal{R}^\omega$ there exists a continuous conjugacy $\eta_{\sigma,\omega}\colon 
  \overline{U_\sigma} \to \tilde \Lambda_\omega$ between
  $f_\sigma$ and $P_\omega$, so that $\eta_{\sigma,\omega}(\mathcal{P}_{\sigma,0}^j)=\mathcal{P}_{0}^j$ for $j\in\{0,...,q-1\}$. The domain $U_\sigma=U_\sigma^n$ for some $n\in \N\cup \{0\}$ and $\overline{U_\sigma}$ contains both critical values $v_1$ and $v_2$. The conjugacy $\eta_{\sigma,\omega}$ is holomorphic in $U_\sigma$.
\end{lemma}

\noindent{\sc Proof: }
Let $f_\sigma \in \mathcal{R}^\omega$. Recall that $\phi_\omega$ and $\phi_\sigma$ are Fatou coordinates for $P_\omega$ and $f_\sigma$ respectively. The map 
$\eta_{\sigma,\omega}=\phi_\omega^{-1}\circ\phi_\sigma \colon
\overline{U_\sigma^0}\to \overline{\mathcal{U}_\omega^0}$, constructed so that $\eta_{\sigma,\omega}(\mathcal{P}_{\sigma,0}^j)=\mathcal{P}_{0}^j$ for all $j\in\{0,...,q-1\}$, is a
homeomorphism, conformal in $U_\sigma^0$ and it
conjugates $f_\sigma$ to $P_\omega$. If $v_2\in\overline{U_\sigma^0}$ then $U_\sigma=U_\sigma^0$ and the proof is done. If not, there exists
$N> 0$ so that $v_2\in \overline{U_\sigma^N}\setminus \overline{U_\sigma^{N-1}}$ and the conjugacy extends, by iterated lifting with respect to the
dynamics, to a conjugacy $\eta_{\sigma,\omega}\colon
\overline{U_\sigma^N}\to \overline{\mathcal{U}_\omega^N}$. Each lift is chosen to agree with the previous map on their common domain of definition.   \qed

\begin{lemma}\label{lemma:obstruction}
  For all non-degenerate parabolic $f_\sigma\in \mathcal{R}^\omega$, 
  $\eta_{\sigma,\omega}(v_2)\in \tilde \Lambda_\omega\setminus S_p$.
\end{lemma}
\noindent{\sc Sketch of Proof: }
The proof is by contradiction. Assume $\exists f_\sigma\in
\mathcal{R}^\omega$ so that $\eta_{\sigma,\omega}(v_2)\in \tilde
\Lambda_\omega\cap S_p$. Let $\overline{U}=\overline{U^n_\sigma}$ be the maximal domain of the
conjugacy $\eta_{\sigma,\omega}$, so that $v_1, v_2 \in \overline{U}$,
and $V=\CC\setminus \overline {U}\cong \D$. Hence $f_\sigma$ has two
univalent inverse branches $f^{-1}_\sigma\colon V\to V$. Let $\mathcal{P}$
denote the union of $q$ repelling petals at the parabolic fixed point
$z_0$, sufficiently small so that $\mathcal{P}\subset
f^{-1}_\sigma(V)\subset V$. If $f_\sigma\in \mathcal{R}^\omega$
so that $\eta_{\sigma,\omega}(v_2)\in \tilde \Lambda_\omega\cap S_p$,
then $f_\sigma^{-1}(\overline U)$ separates $\alpha$ from its
co-preimage $\alpha'$, and $\mathcal{P}$ is then contained in the
image of one of the inverse branches $f^{-1}_\sigma\colon V\to V$. But then
this inverse branch has an attracting $q$-cycle on the ideal
boundary, contradicting the Denjoy--Wolff theorem.
\qed
\vspace{0.5cm}

\noindent{\sc Strategy of Proof Thm. \ref{thm:map}: }
Definition of the map $\chi^{\omega}$. Let $f_\sigma \in
\mathcal{R}^{\omega}$, with second critical value $v_2$. If $f_\sigma$ is degenerate parabolic, then it has two $q$-cycles of components as immediate basin,
with each cycle containing a critical point. In this case choose one of the critical points to be the closest critical point $c_1$, and name the
components in the corresponding cycle $B_0, ..., B_{q-1}$ counter--clockwise, so that
$B_0$ contains the critical point $c_1$. The components of the other
$q$-cycle will be named counter--clockwise so that the component
$B_j'$ is the component between $B_j$ and $B_{(j+1) \mod q}$.  The map $\chi^{\omega}\colon \mathcal{R}^{\omega}\to \hat \Lambda^{\omega}$
is defined by:
\begin{equation}
\chi^{\omega}(\sigma) =  \left\{\begin{array}{ll}
\pi_{\omega}\circ
  \eta_{\sigma,\omega}(v_2)& \text{ if }f_\sigma \text{ non--degenerate}\\
 \pi_{\omega}(\hat\alpha_j)& \text{ if } f_\sigma \text{ degenerate and
 }v_2\in B_j'.
\end{array}\right. 
\end{equation}
The map is well defined by Lemmas \ref{lemma:conjugacy} and
\ref{lemma:obstruction}, and by the equivalence relation
$\sim p/q$, which identifies the images
$\eta_{\sigma,\omega}(v_2)$ and $\hat\alpha_j$'s respectively, in the
cases where there is an ambiguity in the choice of $v_2$. That the map
$\chi^{\omega}$ is holomorphic in the interior will follow from
holomorphic dependence of the Fatou coordinate $\phi_\sigma$ on the
parameter $\sigma$.

Injectivity follows by a classical pull-back argument, see for
example \cite{m92} and \cite{pt09}, adapted to the parabolic situation.
Surjectivity is proved by constructing a sequence of polynomial--like
maps $f_n \in \per(\lambda_n)$, with $\lambda_n\in\D^*$,
$\lambda_n\to \omega$ radially, so that the limiting map $f_\sigma\in
\per(\omega)$ has the correct position of the second critical value.
This is done by using results from \cite{gold90} on the escape loci in $\per(\lambda)$ and results on
convergence of polynomial basins, built upon the star--construction
from \cite{pet99}. Continuity of the map $(\chi^{\omega})^{-1}$ on compact subsets of the bubble--tree is proved by using that  the map $\eta_{\sigma,\omega}$ preserves the
combinatorial structure of the bubble--tree.  
\qed

\section*{Acknowledgements}
Part of this work was done with the support of the ANR grant ``At the
Boundary of Chaos'', Institut de Math\'ematiques de Toulouse. The author would also like to thank Carsten Lunde
Petersen for many productive discussions and Pascale Roesch and Arnaud
Ch\'eritat for helpful comments to the present text.

\end{document}

%% file: figure1.pstex_t
\begin{picture}(0,0)%
\includegraphics{figure1.pstex}%
\end{picture}%
\setlength{\unitlength}{1657sp}%
\begingroup\makeatletter\ifx\SetFigFontNFSS\undefined%
\gdef\SetFigFontNFSS#1#2#3#4#5{%
  \reset@font\fontsize{#1}{#2pt}%
  \fontfamily{#3}\fontseries{#4}\fontshape{#5}%
  \selectfont}%
\fi\endgroup%
\begin{picture}(12159,2991)(511,-3628)
\put(9076,-1696){\makebox(0,0)[lb]{\smash{{\SetFigFontNFSS{5}{6.0}{\rmdefault}{\mddefault}{\updefault}{\color[rgb]{0,0,0}$\pi_\omega$}%
}}}}
\put(10206,-1429){\makebox(0,0)[lb]{\smash{{\SetFigFontNFSS{5}{6.0}{\rmdefault}{\mddefault}{\updefault}{\color[rgb]{0,0,0}$B_2$}%
}}}}
\put(11328,-1768){\makebox(0,0)[lb]{\smash{{\SetFigFontNFSS{5}{6.0}{\rmdefault}{\mddefault}{\updefault}{\color[rgb]{0,0,0}$B_1$}%
}}}}
\put(10986,-2316){\makebox(0,0)[lb]{\smash{{\SetFigFontNFSS{5}{6.0}{\rmdefault}{\mddefault}{\updefault}{\color[rgb]{0,0,0}$B_3$}%
}}}}
\put(12132,-2166){\makebox(0,0)[lb]{\smash{{\SetFigFontNFSS{5}{6.0}{\rmdefault}{\mddefault}{\updefault}{\color[rgb]{0,0,0}$B_0$}%
}}}}
\put(11876,-2731){\makebox(0,0)[lb]{\smash{{\SetFigFontNFSS{5}{6.0}{\rmdefault}{\mddefault}{\updefault}{\color[rgb]{0,0,0}$B_4$}%
}}}}
\put(11502,-908){\makebox(0,0)[lb]{\smash{{\SetFigFontNFSS{5}{6.0}{\rmdefault}{\mddefault}{\updefault}{\color[rgb]{0,0,0}$\mathcal{X}^{\omega_{2/5}}$}%
}}}}
\put(526,-1104){\makebox(0,0)[lb]{\smash{{\SetFigFontNFSS{5}{6.0}{\rmdefault}{\mddefault}{\updefault}{\color[rgb]{0,0,0}$X^{\omega_{2/5}}$}%
}}}}
\put(3203,-1479){\makebox(0,0)[lb]{\smash{{\SetFigFontNFSS{5}{6.0}{\rmdefault}{\mddefault}{\updefault}{\color[rgb]{0,0,0}$B_1$}%
}}}}
\put(1807,-2002){\makebox(0,0)[lb]{\smash{{\SetFigFontNFSS{5}{6.0}{\rmdefault}{\mddefault}{\updefault}{\color[rgb]{0,0,0}$B_3$}%
}}}}
\put(2412,-2823){\makebox(0,0)[lb]{\smash{{\SetFigFontNFSS{5}{6.0}{\rmdefault}{\mddefault}{\updefault}{\color[rgb]{0,0,0}$B_4$}%
}}}}
\put(3756,-2858){\makebox(0,0)[lb]{\smash{{\SetFigFontNFSS{5}{6.0}{\rmdefault}{\mddefault}{\updefault}{\color[rgb]{0,0,0}$B_0$}%
}}}}
\put(2516,-1536){\makebox(0,0)[lb]{\smash{{\SetFigFontNFSS{5}{6.0}{\rmdefault}{\mddefault}{\updefault}{\color[rgb]{0,0,0}$B_2$}%
}}}}
\put(2375,-928){\makebox(0,0)[lb]{\smash{{\SetFigFontNFSS{5}{6.0}{\rmdefault}{\mddefault}{\updefault}{\color[rgb]{0,0,0}$S_p$}%
}}}}
\put(7627,-1532){\makebox(0,0)[lb]{\smash{{\SetFigFontNFSS{5}{6.0}{\rmdefault}{\mddefault}{\updefault}{\color[rgb]{0,0,0}$\hat\alpha_0$}%
}}}}
\put(5757,-1021){\makebox(0,0)[lb]{\smash{{\SetFigFontNFSS{5}{6.0}{\rmdefault}{\mddefault}{\updefault}{\color[rgb]{0,0,0}$\hat{X}^{\omega_{2/5}}$}%
}}}}
\put(8001,-2012){\makebox(0,0)[lb]{\smash{{\SetFigFontNFSS{5}{6.0}{\rmdefault}{\mddefault}{\updefault}{\color[rgb]{0,0,0}$B_0$}%
}}}}
\put(5964,-1723){\makebox(0,0)[lb]{\smash{{\SetFigFontNFSS{5}{6.0}{\rmdefault}{\mddefault}{\updefault}{\color[rgb]{0,0,0}$B_2$}%
}}}}
\put(6191,-2937){\makebox(0,0)[lb]{\smash{{\SetFigFontNFSS{5}{6.0}{\rmdefault}{\mddefault}{\updefault}{\color[rgb]{0,0,0}$B_3$}%
}}}}
\put(7866,-2640){\makebox(0,0)[lb]{\smash{{\SetFigFontNFSS{5}{6.0}{\rmdefault}{\mddefault}{\updefault}{\color[rgb]{0,0,0}$\hat\alpha_4$}%
}}}}
\put(5918,-2350){\makebox(0,0)[lb]{\smash{{\SetFigFontNFSS{5}{6.0}{\rmdefault}{\mddefault}{\updefault}{\color[rgb]{0,0,0}$\hat\alpha_2$}%
}}}}
\put(7025,-1151){\makebox(0,0)[lb]{\smash{{\SetFigFontNFSS{5}{6.0}{\rmdefault}{\mddefault}{\updefault}{\color[rgb]{0,0,0}$B_1$}%
}}}}
\put(7350,-3031){\makebox(0,0)[lb]{\smash{{\SetFigFontNFSS{5}{6.0}{\rmdefault}{\mddefault}{\updefault}{\color[rgb]{0,0,0}$B_4$}%
}}}}
\put(6766,-3124){\makebox(0,0)[lb]{\smash{{\SetFigFontNFSS{5}{6.0}{\rmdefault}{\mddefault}{\updefault}{\color[rgb]{0,0,0}$\hat\alpha_3$}%
}}}}
\put(6401,-1348){\makebox(0,0)[lb]{\smash{{\SetFigFontNFSS{5}{6.0}{\rmdefault}{\mddefault}{\updefault}{\color[rgb]{0,0,0}$\hat\alpha_1$}%
}}}}
\end{picture}%